\documentclass[10pt,twoside]{amsart}
\usepackage{lmodern}
\usepackage[T1]{fontenc}
\usepackage{amsmath,amsfonts,amsthm,amssymb,amscd}
\usepackage{mathrsfs}
\usepackage{latexsym}
\usepackage[pdftex]{graphicx}
\usepackage{enumitem}
\usepackage{bbm}
\usepackage{accents}
\usepackage{euscript}
\usepackage[utf8]{inputenc}

\usepackage{geometry}
\newtheorem{lemma}{Lemma}[section]
\newtheorem{theorem}{Theorem}[section]

\usepackage{hyperref}  
\hypersetup{colorlinks=true, linktoc=all,
    linkcolor=blue,
		citecolor=blue,
		urlcolor=cyan,
		bookmarksopen=true
		}
\usepackage{hyperref}  
\newcommand{\lnc}{\mathscr{L}}
\newcommand{\V}{\EuScript{V}}

\def\M{{\mathfrak M}}
\def\m{{\mathfrak m}}

\def\A{{\mathscr A}}

\def\F{ \mathbb{F}}
\def\R{ \mathbb{R}}
\def\S{\EuScript{S}}

\title{Generalized Estermann problem for non-integer powers with almost proportional summands}
\author{Firuz Rakhmonov}
\address{A.Dzhuraev Institute of Mathematics,  National Academy of Sciences of Tajikistan}
\email{rakhmonov.firuz@gmail.com}
\date{}

\author{Parviz Rakhmonov}
\address{Marex Group, London, UK}
\email{parviz.msu@gmail.com}
\date{}

\begin{document}

\begin{abstract}
For $H \ge N^{1-\frac{1}{2c}} \lnc^{2}$, where $\lnc = \ln N$ and $c$ is a fixed non-integer number satisfying
$$
\|c\| \ge 3c\left(2^{[c]+1}-1\right)\frac{\ln \lnc}{\lnc},
\qquad
c > \frac{4}{3}\left(1 + \frac{52\ln \lnc}{\lnc}\right),
$$
we obtain an asymptotic formula for the number of representations of a sufficiently large integer $N$ in the form
$$
p_{1} + p_{2} + [n^{c}] = N,
$$
where $p_{1}, p_{2}$ are prime numbers, $n$ is a natural number, and
$$
|p_{k} - \mu_{k}N| \le H,\qquad k = 1,2,\qquad |[n^{c}] - \mu_{3}N| \le H,
$$
with $\mu_{1}, \mu_{2}, \mu_{3}$ being fixed positive constants satisfying $\mu_{1} + \mu_{2} + \mu_{3} = 1$. 


\emph{Bibliography:} 21 references.
\end{abstract}

\maketitle

\section{Introduction}
Estermann~\cite{Estermann}, in the case $n=2$, proved an asymptotic formula for the number of solutions of the equation
\begin{equation}\label{formula uravn Estermanna}
 p_1 + p_2 + m^n = N,
\end{equation}
where $p_1$, $p_2$ are prime numbers and $m$ is a natural number. In the works~\cite{RakhmonovZKh-Matzametki-2003-74-4,RakhmonovZKh-Matzametki-2014-95-3,RAO+RFZ-Iss..Saratov-2016-8}, for $n=2,3,4$, this problem was studied under more restrictive conditions, namely, when the summands are almost equal. More precisely, an asymptotic formula was obtained for the number of solutions of the Diophantine equation~(\ref{formula uravn Estermanna}) under the conditions
\[
\left| p_i - \frac{N}{3} \right| \le H,\quad i=1,2,\qquad
\left| m^n - \frac{N}{3} \right| \le H,\qquad H \ge N^{1-\theta(n)} \lnc^{c_n},
\]
respectively, for
\begin{align}\label{formula pannie rezul PrEstcPRS}
\theta(2)=\frac14,\quad c_2=2;\qquad
\theta(3)=\frac16,\quad c_3=3;\qquad
\theta(4)=\frac1{12},\quad c_4=\frac{40}{3}.
\end{align}

In~\cite{RakhmonovFZ-ChebSbornik-2024-25-4}, for an arbitrary fixed $n \ge 3$, based on new estimates for short exponential sums of H.~Weyl~\cite{RakhmonovZKh+FZ-ChebSbornik-2024-25-2,RZKh+RFZ-DNAT-2024-67-3-4,RZKh+RFZ-DNAT-2023-66-9-10,RZKh+RFZ-DNAT-2023-66-11-12,RakhmonovZKh-ChebSbornik-2023-24-3}, the Estermann problem with almost proportional summands was investigated, and an asymptotic formula was obtained for the number of solutions of equation~(\ref{formula uravn Estermanna}) under the conditions
\begin{equation}\label{formula RFZ-ObPrEstSPS}
\left| p_k - \mu_k N \right| \le H,\qquad k=1,2,\qquad
\left| m^n - \mu_3 N \right| \le H,\qquad
H \ge N^{\,1-\frac1{n(n-1)}} \lnc^{\,\frac{2^{n+1}}{n-1} + n - 1}.
\end{equation}

It should be noted that the obtained asymptotic formula, in the case $\mu_1=\mu_2=\mu_3=\tfrac13$, reduces to an asymptotic formula for the generalized Estermann problem with almost equal summands, while the results~\cite{RakhmonovZKh-Matzametki-2014-95-3,RAO+RFZ-Iss..Saratov-2016-8} presented in~(\ref{formula pannie rezul PrEstcPRS}) are particular cases of~(\ref{formula RFZ-ObPrEstSPS}).

V.~N.~Chubarikov posed the following problem: \emph{for fixed non-integer values of $c$, to investigate short exponential sums of the form}
$$
S_c(\alpha; x, y) = \sum_{x-y < n \le x} e\bigl(\alpha [n^{c}]\bigr),
$$
\emph{as well as equation~(\ref{formula uravn Estermanna}), in which the term $m^{n}$ is replaced by $[m^{c}]$, with the case of almost equal summands being considered.}


The first part of this problem was solved in~\cite{RakhmonovPZ-Vestnik MGU 2012,RakhmonovPZDoklAN,RakhmonovPZ-MZ-2014}.
For the sums $S_c(\alpha; x, y)$, an estimate uniform in $c$ was obtained for all $\alpha \in [-0.5, 0.5]$, except for a small neighborhood of zero, and an asymptotic formula with a remainder term, uniform in $c$ for $\alpha$ from a small neighborhood of zero, was also established.
These results made it possible to prove an asymptotic formula in the generalized ternary Estermann problem for non-integer powers with almost equal summands~\cite{RakhmonovPZ-MZ-2016,RakhmonovPZ-ChebSbornik-2015-16-1,RakhmonovPZ-IzvANRT-2013-2}, namely, for $H \ge N^{1-\frac1{2c}} \lnc^2$, an asymptotic formula was obtained for the number of solutions of the equation
$$
p_1 + p_2 + [n^{c}] = N, \qquad
\left| p_i - \frac{N}{3} \right| \le H,\; i=1,2, \qquad
\left| [n^{c}] - \frac{N}{3} \right| \le H,
$$
in prime numbers $p_1$, $p_2$ and a natural number $n$, where $c$ is a fixed non-integer satisfying the conditions
$$
\|c\| \ge 3c\bigl(2^{[c]+1} - 1\bigr)\frac{\ln \lnc}{\lnc},
\qquad
c > \frac{4}{3} + \lnc^{-0.3}.
$$

In the present work, for an arbitrary fixed non-integer $c$ satisfying condition~(\ref{formula usl dlya c}), an asymptotic formula is established in the generalized Estermann problem with almost proportional summands. The obtained result constitutes a generalization and refinement of the main theorem from~\cite{RakhmonovPZ-MZ-2016}.

{\theorem 
Let $N$ be a sufficiently large natural number, $\lnc = \ln N$, and let $\mu_1$, $\mu_2$, $\mu_3$ be positive fixed numbers satisfying $\mu_1 + \mu_2 + \mu_3 = 1$. Let $c$ be a fixed non-integer satisfying the conditions
\begin{equation}\label{formula usl dlya c}
\|c\| \ge 3c\bigl(2^{[c]+1} - 1\bigr)\frac{\ln \lnc}{\lnc}, \qquad
c > \frac{4}{3}\left(1 + \frac{52\ln \lnc}{\lnc}\right).
\end{equation}
Let $J_c(N,H)$ denote the number of solutions of the Diophantine equation
$$
p_1 + p_2 + [n^c] = N
$$
in prime numbers $p_1$, $p_2$ and a natural number $n$ under the conditions
$$
\left| p_k - \mu_k N \right| \le H,\quad k = 1,2, \qquad
\left| [n^c] - \mu_3 N \right| \le H.
$$
Then, for $H \ge N^{\,1-\frac1{2c}} \lnc^{2}$, the following asymptotic formula holds:
$$
J_c(N,H) =
\frac{3H^2}{c(\mu_3 N)^{\,1-\frac1c}\lnc^{2}}
+ O\!\left(\frac{H^2}{N^{\,1-\frac1c}\lnc^{3}}\right),
$$
where the constant implied by the $O$-symbol depends on $\mu_1$, $\mu_2$, $\mu_3$, and $c$.}
\section{Auxiliary Lemmas}
{\lemma \label{Lemma razlozh psi po nuyam zeta(s)} {\rm \cite{Karacuba-OATCh}.} 
Let $2 \le T_0 \le x$, and let $\rho = \beta + i\gamma$ be the nontrivial zeros of the zeta-function. Then
$$
\psi(x)=x-\sum_{|\gamma |\le T_0}\frac{x^{\rho}}{\rho}+R_1(x,T_0),\qquad R_1(x,T_0,)\ll \frac{x\lnc_x^2}{T_0}.
$$}

{\lemma \label{Lemma granitsa nuley zeta(s)} {\rm \cite{Karacuba-OATCh}.} 
There exists an absolute constant $c>0$ such that $\zeta(s)\neq 0$ in the region
$$
\sigma\ge 1-\delta(t),\quad \delta(t)=\frac{c}{\ln^\frac23(2t+2)\ln\ln(2t+2)},
$$}

{\lemma \label{Lemma plotnostn-teor Zhan Tao} {\rm \cite{Zhan Tao-Acta-Math-Sinica-1992}.} 
Let $\varepsilon$ be an arbitrarily small positive constant, and let $T^{\frac{35}{108}+\varepsilon} \le H \le T$. Then the following estimates hold:
$$
N(u,T+H)-N(u,T)\ll\left\{
\begin{array}{ll}
(qH)^{\frac{4}{3-2u}(1-u)}(\ln qH)^9, & \hbox{for} \ \dfrac{1}{2}\le u\le\dfrac{3}{4}, \vspace{7pt} \\
(qH)^{\frac{2}{u}(1-u)+\varepsilon}, & \hbox{for}\ \dfrac{3}{4}\le u\le 1,
\end{array}
\right.
 $$}
{\lemma \label{Lemma Baker i Harman pi(x)-pi(x-y)} {\rm \cite{Baker+Harman1996}.}
Let $y \ge x^{0.534}$. Then the estimate
$$
\frac{y}{\ln x}\ll \pi(x)-\pi(x-y)\ll \frac{y}{\ln x}.
$$}
{\lemma \label{Lemma otsenka S_c(alpha ;x,y)} {\rm \cite{RakhmonovPZ-MZ-2014}.}
Let $x \ge x_0 > 0$, $\lnc_x = \ln x$, and let $A$ be a fixed positive constant greater than one. Let $c$ be a non-integer satisfying the conditions
\begin{equation}\label{uslovie dlya c t vv}
1<c\le \log_2\lnc_x -\log_2 \ln \lnc_x^{6A}, \qquad  \|c\|\ge \left(2^{[c]+1}-1\right)(A+1)\, \lnc_x^{-1} \ln \lnc_x.
\end{equation}
Then, for $y \ge \sqrt{2c x}\, \lnc_x^{A+\theta}$ and $x^{1-c}y^{-1}\lnc_x^A \le |\alpha| \le 0.5$, the estimate
$$
S_c(\alpha ;x,y) \ll y \lnc_x^{-A},
$$
holds, where $\theta = 0$ for $c \ge 1.1$ and $\theta = 0.5$ for $c < 1.1$.}
{\lemma\label{Lemma poveda S_c(alpha ;x,y) v okr nulya} {\rm \cite{RakhmonovPZ-MZ-2014}.}
Let $x \ge x_0 > 0$, and let $A$ be a fixed positive constant greater than one. Let $c$ be a non-integer satisfying conditions~(\ref{uslovie dlya c t vv}). Then, for $y \ge \sqrt{2c}\,x^{\frac{1}{2}}\lnc_x^A$ and $|\alpha| \le x^{1-c}y^{-1}\lnc_x^A$, the following asymptotic formula holds:
$$
S_c(\alpha ;x,y)=\frac{\sin \pi \alpha }{\pi \alpha}\int_{x-y}^xe(\alpha (t^c-0,5))dt +O\left(\frac{y|\sin \pi \alpha |}{\lnc_x^A }\right).
 $$}
{\lemma\label{Lemma S(alpha;Nk,H=S1(alpha;Nk+H,2H)} {\rm \cite{RakhmonovFZ-ChebSbornik-2024-25-4}}
Let $\mu_k$ be a fixed real number with $0<\mu_k<1$, let $N$ be a sufficiently large natural number, and let $N_k=\mu_k N+H$, $k=1,2$, with $N^{\frac12}\le H\le N^{1-\frac1{30}}$. Define
$$
\S(\alpha;N_k,2H)=\sum_{N_k-2H<p\le N_k}e(\alpha p),\qquad S_1(\alpha;x,y)=\sum_{x-y<n\le x}\Lambda(n)e(\alpha n).
$$
Then the relation
$$
\S(\alpha;N_k,2H)=\frac{S_1(\alpha;N_k,2H)}{\ln(\mu_kN)}+O\left(\frac{H^2}{N\ln(\mu_kN)}\right).
$$
holds.}
{\lemma \label{Lemma poved LinKorTrigSum v okr nulya}
Let $x\ge x_0$, let $A$ be an arbitrary fixed positive constant, let $y\ge x^{\frac{5}{8}}\lnc_x^{{1.5}A+15} x$, and let $|\alpha| \le x(2 \pi y^2)^{-1}$. Then the following asymptotic formula holds:
$$
S_1(\alpha;x,y)=\sum_{x-y<n\leq x}\Lambda(n)e(\alpha n)=\frac{\sin\pi\alpha y}{\pi\alpha}e\left(\alpha\left(x-\frac{y}{2}\right)\right)+
O\left(\frac{y}{\lnc_x^A}\right).
$$}
{\sc Proof.} Without loss of generality, we may assume that
\begin{equation}\label{formula y=x58..}
y=x^{\frac58}\lnc_x^{1.5A+15}.
\end{equation}
Applying Abel's transformation in integral form, we obtain:
\begin{align*}
S_1(\alpha;x,y)
&=-\int_{x-y}^{x}\psi (u)de(\alpha u)+e(\alpha x)\psi (x)- e(\alpha (x-y))\psi (x-y ).
\end{align*}
Using the representation of the Chebyshev function as a sum over the zeros of the zeta-function (Lemma~\ref{Lemma razlozh psi po nuyam zeta(s)}) with $T_0=\left(xy^{-1}+|\alpha|x\right)\lnc_x^{A+2}$, we find:
\begin{align*}
S_1(\alpha;x,y)&=-\int_{x-y}^x\left(u-\sum_{|\gamma |\le T_0}\frac{u^\rho}{\rho}\right)de(\alpha u)+e(\alpha x)\left(x-\sum_{|\gamma|\le T_0}\frac{x^\rho}{\rho}\right)-\\
&-e(\alpha(x-y))\left((x-y)-\sum_{|\gamma|\le T_0}\frac{(x-y)^\rho}{\rho}\right)-\\
&-\int_{x-y}^xR(u,T_0)2\pi i\alpha e(\alpha u)du+e(\lambda x)R(x,T_0)-e(\alpha(x-y))R(x-y,T_0).
\end{align*}
Applying integration by parts to the first integral, and also using the estimate for $R(u,T_0)$ from Lemma~\ref{Lemma razlozh psi po nuyam zeta(s)} together with the value of the parameter $T_0$, we obtain
\begin{align}
S_1(\alpha;x,y)&=\int_{x-y}^xe(\lambda u)du-\sum_{|\gamma|\le T_0}\int_{x-y}^xu^{\rho-1}e(\lambda u)du+ O\left(\frac{y}{q^\frac12\lnc_x^A}\right)=\nonumber\\
&=\frac{\sin\pi\lambda y}{\pi\lambda}e\left(\lambda\left(x-\frac{y}{2}\right)\right)- W(\alpha,x,y)+ O\left(\frac{y}{\lnc_x^A}\right),\label{formula S(alpha;x,y)=..}
\end{align}
where
$$
W(\alpha,x,y)=\sum_{|\gamma |\le T_0}I(\rho),\qquad I(\rho)=\int_{x-y}^xu^{\beta-1}e\left(\alpha u+\frac\gamma{2\pi}\ln u\right)du.
$$
The sum $W(\alpha,x,y)$ will be estimated only in the case $\alpha \ge 0$. For $\alpha \le 0$, using the relation
\begin{align*}
\overline{W(\alpha,x,y)}&= \sum_{|\gamma |\le T_0}\int_{x-y}^{x}u^{\beta-1}e\left(-\alpha u-\frac\gamma{2\pi}\ln u\right)du=\\
&=\sum_{|\overline{\gamma} |\le T_0}\ \int_{x-y}^{x}u^{\rho -1}u^{\beta-1}e\left(-\alpha u+\frac\gamma{2\pi}\ln u\right)du=W(-\alpha,x,y),
\end{align*}
we reduce its estimation to the case $\alpha\ge 0$. Estimating the integral $I(\rho)$ trivially, as well as via the magnitude of the first derivative (see~\cite{VoroninKaratsuba}, p.~359), we obtain:
\begin{equation}\label{Otsenka int I(rho)}
|I(\rho )| \ll x^\beta \min_{x-y\le u\le x}\left( \frac{y}{x},\frac{1}{\min |\gamma+2\pi \alpha u|} \right).
\end{equation}
All zeros $\rho =\beta+i\gamma$ satisfying $|\gamma|\le T_0 $ are partitioned into the sets $D_1$, $D_2$, and $D_3$ as follows:
\begin{align*}
&D_1=\left\{\rho:  -T_0 +2\pi\alpha u \le\gamma+2\pi\alpha u<-2\pi\alpha x+2\pi\alpha u-\frac{x}{y}\right\}, \\
&D_2=\left\{\rho:-2\pi\alpha x+2\pi\alpha u-\frac{x}{y}\le\gamma+2\pi\alpha u\le2\pi\alpha u-2\pi\alpha(x-y)+\frac{x}{y}\right\} ,\\
& D_3=\left\{\rho:\ 2\pi\alpha u-2\pi\alpha(x-y)+\frac{x}{y}<\gamma+2\pi\alpha u\le T_0+2\pi\alpha u\right\}.
\end{align*}
Denoting by $W_1$, $W_2$, and $W_3$ the sums of the absolute values of the integral $I(\rho )$ over zeros belonging to the sets $D_1$, $D_2$, and $D_3$, respectively, we have:
\begin{equation}
\label{Formula W(alpha,x,y)<<W1+W2+W3.}
|W(\alpha,x,y)|\le \sum_{|\gamma |\le T_0}|I(\rho)|=W_1+W_2+W_3.
\end{equation}
On the interval $x-y\leq u \le x$, the function $2\pi \alpha u$ is monotonically increasing; therefore, for the right boundary of $D_1$ and the left boundary of $D_3$, respectively, we have
$$
-2\pi\alpha x+2\pi\alpha u-\frac{x}{y}\le-\frac{x}{y},\quad 2\pi\alpha u-2\pi\alpha(x-y)+\frac{x}{y}\ge\frac{x}{y}.
$$
Consequently, if $\rho$ belongs to $D_1$, $D_2$, or $D_3$, respectively, then
$$
\gamma+2\pi\alpha u<-\frac{x}{y}, \qquad -\frac{x}{y}\le\gamma+2\pi\alpha u\le\frac{x}{y}, \qquad \gamma+2\pi\alpha u>\frac{x}{y}.
$$
Therefore, for the monotonically increasing function $\gamma +2\pi \alpha u$ on the interval $x-y\le u\le x$, the following relations hold:
\begin{align*}
\min_{x-y\le u\le x}|\gamma+2\pi\alpha u|=-\max_{x-y\le u\le x}(\gamma+2\pi\alpha u)=-\gamma-2\pi\alpha x&\ge\frac{x}{y},\quad \text{if }  \quad \rho\in D_1, \\
-\frac{x}{y}-2\pi\alpha y\le \gamma +2\pi \alpha u  \le \frac{x}{y}+&2\pi\alpha y, \quad \text{if }  \quad \rho\in D_2, \\
\min_{x-y\le u\le x}|\gamma+2\pi\alpha u|=\min_{x-y\le u\le x}(\gamma+2\pi\alpha u)=\gamma+2\pi\alpha(x-y)&\ge\frac{x}{y},\quad \text{if }  \quad \rho\in D_3.
\end{align*}
Hence, taking into account estimate~(\ref{Otsenka int I(rho)}), for $W_1$, $W_2$, and $W_3$ we obtain
$$
W_1\ll\sum_{\rho\in D_1}\frac{x^\beta}{-\gamma-2\pi\alpha x}, \qquad W_2\ll\sum_{\rho\in D_2}yx^{\beta-1}, \qquad W_3\ll\sum_{\rho \in D_3}\frac{x^\beta}{\gamma+2\pi\alpha(x-y)}.
$$
The sums $W_1$ and $W_3$ are estimated in the same way. Let us estimate $W_1$. All zeros in the set
$$
D_1=\left\{\rho:\ \frac{x}{y}<-\gamma-2\pi\alpha x\le T_0-2\pi\alpha x\right\},
$$
are divided into classes $D_{11},\ldots,D_{1r}$, $r\ll \ln T_0\ll \ln x$, as follows: the class $D_{1n}$ consists of those zeros $\rho$ for which the conditions
$$
\frac{nx}y<-\gamma-2\pi\alpha x\le\frac{(n+1)x}y.
$$
hold. Therefore,
\begin{align*}
W_1&\ll\sum_{n=1}^r\sum_{\rho\in D_{1n}}\frac{x^\beta}{-\gamma-2\pi\lambda x} \le\frac{y}{x}\sum_{n=1}^r\sum_{\rho\in D_{1n}}\frac{x^\beta }{n}\le \frac{y\lnc_x}x\max_{1\le n\le r}\sum_{\rho\in D_{1n}}x^{\beta}\le\frac{y\lnc_x}x\max_{|T|\le T_0}\sum_{T-\frac{x}y<\gamma\le T}x^\beta.
\end{align*}
Let us estimate $W_2$. Representing the set $D_2$ in the form
\begin{align*}
D_2
&=\left\{\rho:\ T_1-2\pi\alpha y-\frac{2x}{y}\le-\gamma\le T_1\right\}, \qquad T_1=2\pi\alpha x+\frac{x}{y}\le T_0,
\end{align*}
and taking into account that for $\alpha\le x\left(2\pi y^2\right)^{-1}$ the length of the set $D_2$ satisfies the inequality
$$
2\pi\lambda y+\frac{2x}{y}\le\frac{3x}{y},
$$
and also using the trivial estimate for the integral $I(\rho,\lambda)$, i.e., the first estimate in~(\ref{Otsenka int I(rho)}), we obtain
\begin{align*}
W_2&\le\sum_{\rho \in D_2}|I(\rho)|\ll\frac{y}{x}\sum_{\rho \in D_2}x^\beta\le\frac{3y}{x}\max_{|T|\le T_0}\sum_{T-\frac{x}{y}\le-\gamma\le T}x^\beta\ll\frac{y}{x}\max_{|T|\le T_0}\sum_{T-\frac{x}{y}\le\gamma\le T}x^\beta.
\end{align*}
Substituting the obtained estimates for $W_1$, $W_2$, and $W_3$ into~(\ref{Formula W(alpha,x,y)<<W1+W2+W3.}), we obtain:
\begin{equation}
\label{Formula W(alpha,x,y)<ylnx/x max W_4}
W(\alpha,x,y)\ll \frac{y\lnc_x}{x} \max_{|T|\le T_0 }\V, \qquad \V=\sum_{T-\frac{x}{y}\le\gamma\le T}x^{\beta }.
\end{equation}
The estimation of the sum $\V$ is reduced to estimating the number of zeros of the Riemann zeta-function in narrow rectangles of the critical strip. We have
\begin{align*}
\V=&\sum_{T-xy^{-1}<\gamma\le T}\left(\int_0^\beta x^udu+1\right)=\lnc_x\int_0^1x^u\sum_{\substack{T-xy^{-1}<\gamma\le T\\\beta\ge u}}du+\sum_{T-xy^{-1}<\gamma\le T}1\\
&=\lnc_x\int_0^1x^u\left(N(u,T)-N(u,T-xy^{-1})\right)du+\left(N(T)-N(T-xy^{-1})\right).
\end{align*}
Hence, taking into account that the nontrivial zeros $\rho = \beta + i\gamma$ are symmetrically distributed with respect to the critical line $\sigma = 0.5$, we obtain
\begin{align*}
\V&\le\lnc_x\int_{0,5}^1x^u\left(N(u,T)-N(u,T-xy^{-1})\right)du+\left(\frac{\sqrt{x}\lnc_x}{2}+1\right)\left(N(T)-N(T-xy^{-1})\right)\le\\
&\le\lnc_x\max_{u\ge 0,5}x^u\left(N(u,T)-N(u,T-xy^{-1},\chi)\right)+\left(\frac{\sqrt{x}\lnc_x}{2}+1\right)\left(N(T)-N(T-xy^{-1})\right)\le\\
&\le2\lnc_x\max_{u\ge 0,5}x^u\left(N(u,T)-N(u,T-xy^{-1})\right).
\end{align*}
According to Lemma~\ref{Lemma granitsa nuley zeta(s)}, the function $\zeta(\sigma+it)$ has no zeros in the region
$$
\sigma\ge 1-\delta(t),\qquad \delta(t)=\frac{c}{\ln^\frac23(2t+2)\ln\ln(2t+2)}.
$$
Consequently, taking into account that $\delta(T) \ge \delta(T_0)$, we obtain
$$
\V\le2\lnc_x\max_{0,5\le u\le1-\delta}x^u\left(N(u,T)-N(u,T-xy^{-1})\right),\qquad \delta=\delta(T_0).
$$

Substituting the right-hand side of this inequality into~(\ref{Formula W(alpha,x,y)<ylnx/x max W_4}), we obtain
\begin{align}\label{Formula W(alpha,x,y)<ylnx/x(N(u,T)-N(u,t-x/y))}
|W(\alpha;x,y)|\ll\frac{y\lnc_x}{x}\max_{|T|\le T_0}\max_{0,5\le u\le1-\delta}x^u\left(N(u,T)-N\left(u,T-\frac{x}{y},\chi\right)\right).
\end{align}
From the relations $|T|\le T_0=\left(xy^{-1}+|\alpha|x\right)\lnc_x^{A+2}$, $0\le\alpha\le x\left(2\pi y^2\right)^{-1}$, and condition~(\ref{formula y=x58..}), we have
\begin{align*}
\frac{T}{(xy^{-1})^3}\le&\left(\frac{y^2}{x^2}+\alpha\frac{y^3}{x^2}\right)\lnc_x^{A+2}\le\left(\frac{y^2}{x^2}+\frac{y}{2\pi x}\right)\lnc_x^{A+2}\le1.
\end{align*}
Consequently, for the right-hand side of inequality~(\ref{Formula W(alpha,x,y)<ylnx/x(N(u,T)-N(u,t-x/y))}), the condition $\frac{x}{y} \ge T^{\frac13}$ holds; that is, the density theorem for narrow rectangles of the critical strip (Lemma~\ref{Lemma plotnostn-teor Zhan Tao}) can be applied to this sum. Setting in this lemma $\varepsilon = \frac{\delta}{3} - \frac{\delta^2}{1-\delta}$, we obtain
\begin{align}
&|W(\alpha;x,y)| \ll \A_1 + \A_2, \label{Formula W(alpha;x,y)<<A1+A2l} \\
&\A_1 = \frac{y \lnc_x^{10}}{x} \max_{0.5 \le u \le 0.75} x^u \left( \frac{x}{y} \right)^{\frac{4-4u}{3-2u}}, \nonumber \\
&\A_2 = \frac{y \lnc_x}{x} \max_{0.75 \le u \le 1-\delta} x^u \left( \frac{x}{y} \right)^{\frac{2}{u}(1-u) + \varepsilon}. \nonumber
\end{align}

\textbf{Estimate of $\A_1$.} We have
\begin{align*}
&\A_1=\frac{x\lnc_x^{10}}{y}\max_{0,5\le u\le0,75}f_1(u),\qquad f_1(u)=x^u\left(\frac{x}{y}\right)^\frac1{u-1,5}>0,\\
&f_1'(u)=f_1(u)\left(\ln x+\frac{\ln\left(\frac{y}{x}\right)}{(u-1,5)^2}\right)=
\frac{f_1(u)}{(u-1,5)^2}\ln\frac{y}{x^{-u^2+3u-1,25}}.
\end{align*}
From condition~(\ref{formula y=x58..}) and the relation
$$
\max_{0,5\le u\le0,75}\left(-u^2+3u-1,25\right)=\left.\left(-u^2+3u-1,25\right)\right|_{u=0,75}=\frac{7}{16},
$$
it follows that
$$
\ln\frac{y}{x^{-u^2+3u-1,25}}\ge\ln\frac{yx^\frac58\lnc_x^{1,5A+15}}{x^\frac7{16}}=\ln x^\frac3{16}\lnc_x^{1,5A+15}>\ln x^\frac15>0.
$$
Hence, in turn, we obtain that on the interval $0,5\le u\le 0,75$ the inequality
$$
f_1'(u)\ge\frac{f(u)}{(u-1,5)^2}\ln x^\frac15>0,
$$
holds, that is, $f_1'(u)$ is positive and $f_1(u)$ is an increasing function on the interval $0,5\le u\le0,75$. Using this property, and then relation~(\ref{formula y=x58..}), we obtain
\begin{align*}
\A_1=\frac{x\lnc_x^{10}}{y}f_1(0.75)=\frac{x\lnc_x^{10}}{y}x^\frac34\left(\frac{x}{y}\right)^{-\frac43}= y\cdot\left(\frac{x^\frac58\lnc_x^{15}}{y}\right)^\frac23=y\lnc_x^{-A}.
\end{align*}

\textbf{Estimate of $\A_2$.} We have
\begin{align*}
\A_2=&\frac{y^{3-\varepsilon}\lnc_x}{x^{3-\varepsilon}}\max_{0,75\le u\le1-\delta}f_2(u),\qquad f_2(u)=x^u\left(\frac{x}{y}\right)^\frac2u>0,\qquad \delta=\delta(T_0).\\
&f_2'(u)=f_2(u)\left(\ln x+\frac{\ln\left(\frac{y}{x}\right)}{u^2}\right)=\frac{f_2(u)}{u^2}\ln\frac{y}{x^{1-u^2}}.
\end{align*}
From condition~(\ref{formula y=x58..}) and the relation
$$
\max_{0,75\le u\le1-\delta}\left(-u^2+1\right)=\left.\left(-u^2+1\right)\right|_{u=0,75}=\frac{7}{16},
$$
it follows that
$$
\frac{y}{x^{1-u^2}}\ge\frac{x^{\frac58}\lnc_x^{1,5A+15}}{x^\frac7{16}}\ge x^\frac3{16}\lnc_x^{1,5A+15}.
$$
Hence, in turn, we obtain that on the interval $0,75\le u\le 1-\delta$ the inequality
$$
f_2'(u)\ge\frac{f_2(u)}{u^2}\ln x^\frac3{16}>0,
$$
holds, that is, $f_2'(u)$ is positive and $f_2(u)$ is an increasing function on the interval $0,75\le u\le1-\delta$. Using this property, we obtain
\begin{align*}
\A_2=&\frac{y^{3-\varepsilon}\lnc_x}{x^{3-\varepsilon}}x^{1-\delta}\left(\frac{x}{y}\right)^\frac2{1-\delta}= y\cdot\frac{x^{\frac{2\delta}{1-\delta}-\delta+\varepsilon}\lnc_x}{y^{\frac{2\delta}{1-\delta}+\varepsilon}} =y\lnc_x\left(\frac{x^{\frac{\delta+\delta^2+(1-\delta)\varepsilon}{2\delta+(1-\delta)\varepsilon}} }{y}\right)^{\frac{2\delta}{1-\delta}+\varepsilon}\hspace{-23pt}=\\
&=y\lnc_x\left(\frac{x^{\frac{\scriptstyle5}{\scriptstyle8}}\lnc_x^{1,5A+15}}y x^{f(\delta,\varepsilon)}\lnc_x^{-1,5A-15}\right)^{\frac{\scriptstyle{2\delta}}{\scriptstyle{1-\delta}}+\varepsilon},\qquad f(\delta,\varepsilon)=\frac{\delta+\delta^2+(1-\delta)\varepsilon}{2\delta+(1-\delta)\varepsilon}-\frac58.
\end{align*}

From this, and from relation~(\ref{formula y=x58..}), we obtain
\begin{align*}
\A_2&\ll y\cdot x^{f(\delta,\varepsilon)\frac{\scriptstyle{2\delta+(1-\delta)\varepsilon}}{\scriptstyle{1-\delta}}}\lnc_x.
\end{align*}
Further, for $\varepsilon=\frac{\delta}3-\frac{\delta^2}{1-\delta}$, using the relation
\begin{align*}
&f(\delta,\varepsilon)\frac{2\delta+(1-\delta)\varepsilon}{1-\delta}
 =\frac{\delta+\delta^2+(1-\delta)\varepsilon}{1-\delta}- \frac{5(2\delta+(1-\delta)\varepsilon)}{8(1-\delta)}=\\
&=\frac{-2\delta+8\delta^2+3(1-\delta)\varepsilon}{8(1-\delta)}=-\frac{\delta}{8}+\frac{3(1-\delta)\varepsilon-\delta+7\delta^2}{8(1-\delta)}= -\frac{\delta}{8},
\end{align*}
we find
\begin{align}
\A_2&\ll yx^{-0,125\delta}\lnc_x=y\lnc_x\exp\left(-0,125\delta\lnc\right).\label{formula A2<<}
\end{align}
Using the conditions $\alpha\le x\left(2\pi y^2\right)^{-1}$ and~(\ref{formula y=x58..}), we have
\begin{align*}
T_0=&\left(\frac xy+\alpha x\right)\lnc_x^{A+3}\le\left(\frac{x}y+\frac{x^2}{2\pi y^2}\right)\lnc_x^{A+3}\le \frac{x^2}{y^2}\lnc_x^{A+3} =x^{\frac34}\lnc_x^{-2A-33}<x^{\frac45},
\end{align*}
Using this inequality, we estimate the parameter $\delta=\delta(T_0)$ from below:
\begin{align*}
\delta(T_0)&=\frac{c_1}{\left(\ln(2T_0+2)\right)^\frac23\ln\ln(2T_0+2)}\ge\frac{c_1}{\lnc_c^\frac23\ln\lnc_x}\ge c_1\lnc_x^{-0.6}.
\end{align*}
Therefore, from~(\ref{formula A2<<}) we obtain
\begin{align*}
\A_2&\ll y\lnc_x\exp\left(-0,125c_1\lnc_x^{0.4}\right)\ll y\lnc_x^{-A}.
\end{align*}

Substituting this estimate and the estimate for the sum $\A_1$ into (\ref{Formula W(alpha;x,y)<<A1+A2l}), we obtain
\begin{align*}
|W(\alpha;x,y)|&\ll y\lnc_x^{-A}.
\end{align*}
From this estimate and (\ref{formula S(alpha;x,y)=..}) we obtain the assertion of Lemma \ref{Lemma poved LinKorTrigSum v okr nulya}.

\subsection{Proof of the theorem}
Without loss of generality, we assume that
\begin{equation}\label{formula H=N..}
H=N^{1-\frac{1}{2c}}\lnc^2
.
\end{equation}
Introducing the notation $N_3^c=\mu_3N+H$ and  $(N_3-H_3)^c=\mu_3N-H$, and then using the definition of the sum $S_c(\alpha;x,y)$, we obtain
\begin{align}
&\sum_{|[n^c]-\mu_3N|\le H}e(\alpha [n^c])=\sum_{(N_3-H_3)^c+\{n^c\}\le n^c\le N_3^c+\{n^c\}}e(\alpha [n^c])=\nonumber\\
&=\sum_{(N_3-H_3)^c<n^c\le N_3^c}e(\alpha [n^c])+\theta_{31}+\theta_{32}=S_c\left(\alpha;N_3,H_3\right)+\theta_3,
\label{formula ..=Sc(alpha;N3,H3)}
\end{align}
where $\theta_3=\theta_{31}+\theta_{32}$, and $\theta_{31}$ and $\theta_{32}$ are defined as follows:
\begin{itemize}
  \item $|\theta_{31}|=1$, if there exists an integer in the interval $\left[N_3-H_3, \ \left((N_3-H_3)^c+\{n^c\}\right)^\frac1c\right]$, whose length is less than one, and $|\theta_{32}|=0$ otherwise;
  \item $|\theta_{32}|=1$, if there exists an integer in the half-interval $\left(N_3, \ \left(N_3^c+\{n^c\}\right)^\frac1c\right]$, whose length is also less than one, and $|\theta_{32}|=0$ otherwise.
\end{itemize}
For the parameters $N_3$ and $H_3$, the following relations hold, which will be repeatedly used below:
\begin{align}
&N_3=(\mu_3N+H)^\frac1c=(\mu_3N)^\frac1c\left(1+\frac{H}{\mu_3N}\right)^\frac1c= (\mu_3N)^\frac1c\left(1+\frac{H}{c\mu_3N}+O\left(\frac{H^2}{N^2}\right)\right),\label{formula poryadok N3}\\
&H_3=(\mu_3N+H)^\frac1c-(\mu_3N-H)^\frac1c=\frac{2H}{c(\mu_3N)^{1-\frac1c}}+O\left(\frac{H^2}{N^{2-\frac1c}}\right).\label{formula poryadok H3}
\end{align}
Using relation (\ref{formula ..=Sc(alpha;N3,H3)}) and the notation
$$
\S_1(\alpha;N_k,2H)=\sum_{N_k-2H<p\le N_k}e(\alpha p),\qquad N_k=\mu_kN+H,\qquad k=1,\ 2;
$$
$J_c(N,H)$~---~the number of solutions of the Diophantine equation
$$
p_1+p_2+[n^c]=N,
$$
in prime numbers $p_1$, $p_2$ and natural numbers $n$ under the conditions
$$
|p_k-\mu_kN|\le H, \quad k=1,2,\quad |[n^c]-\mu_3N|\le H,
$$
can be represented in the form
\begin{align}
J_c&(N,H)=\int_{-0.5}^{0,5}e(-\alpha N)\sum_{|p_1-\mu_1 N|\le H}e(\alpha p_1)\sum_{|p_2-\mu_2 N|\le H}e(\alpha p_2)\sum_{|[n^c]-\mu_3N|\le H}e(\alpha[n^c])=\nonumber\\
&=\int_{-0.5}^{0.5}e(-\alpha N)\left(\S_1(\alpha;N_1,2H)+\theta_1\right)\left(\S_1(\alpha;N_2,2H)+\theta_2\right) \left(S_c(\alpha;N_3,H_3)+\theta_3\right)d\alpha,\label{formula Jn(N,H)-1}
\end{align}
where $|\theta_k|$, (for  $k=1,\ 2$) is equal to $1$ if, respectively, the lower bounds of the exponential sums $\S_1(\alpha;N_1,2H)$, $\S_1(\alpha;N_2,2H)$, $T(\alpha;N_3,H_3)$, that is, the numbers $N_1-2H$, $N_2-2H$, $N_3-H_3$, are integers, and is equal to $0$ otherwise. Expanding the brackets in the integrand of (\ref{formula Jn(N,H)-1}), we obtain
\begin{align}
J_c(N,H)=\int_{-0.5}^{0.5}e(-\alpha N)&\S_1(\alpha;N_1,2H)\S_1(\alpha;N_2,2H)S_c(\alpha;N_3,H_3)d\alpha+\R_1,\label{formula Jn(N,H)-2}\\
\R_1=\int_{-0.5}^{0.5}e(-\alpha N)&\left(
\theta_3\S_1(\alpha;N_1,2H)\S_1(\alpha;N_2,2H)+ \theta_1\theta_2S_c(\alpha;N_3,H_3)+\right.\nonumber\\
&\left.+\theta_2\S_1(\alpha;N_1,2H)S_c(\alpha;N_3,H_3)+\theta_1\theta_3\S_1(\alpha;N_2,2H)+\right. \nonumber\\
&\left.+\theta_1\S_1(\alpha;N_2,2H)S_c(\alpha;N_3,H_3)+\theta_2\theta_3\S_1(\alpha;N_1,2H)\right)d\alpha.\nonumber
\end{align}

In $\R_1$, passing to estimates, using the Cauchy inequality, and then the relations
\begin{align*}
&\int_{-0.5}^{0.5}|\S_1(\alpha;N_k,2H)|^2d\alpha=\pi(N_k)-\pi(N_k-2H)\le 2H, \qquad k=1,\ 2;\\
&\int_{-0.5}^{0.5}|S_c(\alpha;N_3,H_3)|^2d\alpha=[N_3]-[N_3-H_3]\le H_3+1\ll\frac{H}{N^{1-\frac1c}},
\end{align*}
where in deriving the second of them relation (\ref{formula poryadok H3}) is used, we obtain
\begin{align*}
\R_1&\le\left(\int_{-0.5}^{0.5}|\S_1(\alpha;N_1,2H)|^2d\alpha\int_{-0.5}^{0.5}|\S_1(\alpha;N_2,2H)|^2d\alpha\right)^\frac12+ \left(\int_{-0.5}^{0.5}|S_c(\alpha;N_3,H_3)|^2d\alpha\right)^\frac12+\\
&+\left(\int_{-0.5}^{0.5}|\S_1(\alpha;N_1,2H)|^2d\alpha\int_{-0.5}^{0.5}|S_c(\alpha;N_3,H_3)|^2d\alpha\right)^\frac12+ \left(\int_{-0.5}^{0.5}|\S_1(\alpha;N_2,2H)|^2d\alpha\right)^\frac12+\\
&+\left(\int_{-0.5}^{0.5}|\S_1(\alpha;N_2,2H)|^2d\alpha\int_{-0.5}^{0.5}|S_c(\alpha;N_3,H_3)|^2d\alpha\right)^\frac12+ \left(\int_{-0.5}^{0.5}|\S_1(\alpha;N_1,2H)|^2d\alpha\right)^\frac12\ll\\
&\ll H+\left(\frac{H}{N^{1-\frac1c}}\right)^\frac12+\left(\frac{H^2}{N^{1-\frac1c}}\right)^\frac12+H^\frac12\ll H\ll\frac{H^2}{N^{1-\frac1c}\lnc^3}.
\end{align*}
From this and from (\ref{formula Jn(N,H)-2}), we obtain
\begin{align*}
J_c(N,H)=&\int_{-0.5}^{0.5}\F(\alpha)e(-\alpha N)d\alpha+O\left(\frac{H^2}{N^{1-\frac1c}\lnc^3}\right),\\
\F(\alpha)=&\S_1(\alpha;N_1,2H)\S_1(\alpha;N_2,2H)S_c(\alpha;N_3,H_3).
\end{align*}
We divide the interval of integration $[-0,5,0,5]$ into points of two classes. To the points of the first class we assign the interval
$$
\M=[-\varkappa,\varkappa], \quad \text{  where } \quad \varkappa  =(2cH)^{-1}\,\lnc^2.
$$
The remaining intervals
$$
\m_+=[\varkappa,0,5]\quad \text{  and } \quad \m_-=[-0,5,-\varkappa ]
$$
are assigned to the points of the second class.

Denote by $I(\M)$,  $I(\m_+)$ and $I(\m_-)$ the integrals over the sets $\M$, $\m_+$ and $\m_-$, respectively. Then we have
$$
J_c(N,H)=I(\M)+I(\m_+)+I(\m_-).
$$
In the last formula, the first term, namely $I(\M)$, provides the main term of the asymptotic formula for $J_c(N,H)$, while $I(\m_+)$ and $I(\m_-)$ are included in its remainder term.

\subsection{Estimate of the integrals $I(\m_+ )$ and $I(\m_-)$}
\label{subsec_tri}
We have
$$
I(\m_+)=\int\limits_{\m_+}e(-\alpha N)\S_1(\alpha;N_1,2H)\S_1(\alpha;N_2,2H)S_c(\alpha;N_3,H_3)d\alpha.
$$
Passing to estimates, and then applying the Cauchy inequality for integrals, we find
\begin{align*}
&I(\m_+)\ll\max_{\alpha \in\m_+}|T(\alpha;N_3,H_3)|\int_0^1|\S_1(\alpha;N_1,2H)||\S_1(\alpha;N_2,2H)|d\alpha =\\
&=\max_{\alpha\in\m_+}|S_c(\alpha;N_3,H_3)|\left(\int_0^1|\S_1(\alpha;N_1,2H)|^2d\alpha\right)^\frac12 \left(\int_0^1|\S_1(\alpha;N_2,2H)|^2d\alpha\right)^\frac12 =\\
&=\max_{\alpha\in\m_+}|S_c(\alpha;N_3,H_3)|\left(\pi\left(\mu_1N+H\right)-\pi\left(\mu_1N-H\right)\right)^\frac12 \left(\pi\left(\mu_2N+H\right)-\pi\left(\mu_2N-H\right)\right)^\frac12.
\end{align*}

Applying Lemma \ref{Lemma Baker i Harman pi(x)-pi(x-y)} to the last two factors on the right-hand side of the obtained formula, taking into account the relation
$$
H=N^{1-\frac{1}{2c}}\lnc^2\ge N^\frac58\lnc^2 \ge \left(\mu_kN+H\right)^{0.534},\qquad k=1,
$$
we find
$$
\pi(\mu_kN+H)-\pi(\mu_kN-H)\ll\frac{H}{\lnc }.
$$
Consequently
\begin{equation}\label{Integral I(M2)}
I(\m_+) \ll\frac{H}\lnc\max_{\alpha \in \m_+}|S_c(\alpha;N_3,H_3)|.
\end{equation}
We estimate $S_c(\alpha, N_3, H_3)$ for $\alpha$ from the set $\m_+$ using Lemma \ref{Lemma otsenka S_c(alpha ;x,y)} in the case $c > 1,1$ with
$$
A = 2, \qquad x = N_3, \qquad y = H_3,
$$
and we also verify each of the following conditions:
\begin{align}
& \|c\| \ge 3\left(2^{[c]+1} - 1\right) \frac{\ln \ln N_3}{\ln N_3}, \label{uslovie1} \\
& H_3 \ge \sqrt{2 c N_3} (\ln N_3)^2, \label{uslovie2} \\
& \varkappa = (2 c H)^{-1} \lnc^2 \ge \frac{(\ln N_3)^2}{H_3 N_3^{\,c-1}}. \label{uslovie3}
\end{align}

Using the definition of the parameter $N_3=(\mu_3N+H)^\frac1c$, we have
\begin{align}
\ln N_3=&\frac1c\ln\left(N\left(\mu_3+\frac{H}{N}\right)\right)=\frac1c\left(\lnc+\ln\left(\mu_3+\frac{H}{N}\right)\right)=\nonumber\\ &=\frac{\lnc}c\left(1+\frac{\ln\left(\mu_3+\frac{H}{N}\right)}{\lnc}\right)=\frac{\lnc}c\left(1+\R(n)\right),\qquad \R(n)\ll\lnc^{-1}. \label{formula poryadok lnN3}
\end{align}
Taking logarithms of both sides of the last equality, we obtain
\begin{align}\label{formula poryadok lnlnN3}
\ln\ln N_3=\ln\lnc-\ln c+\ln\left(1+\frac{\ln\left(\mu_3+\frac{H}{N}\right)}{\lnc}\right)=\ln\lnc+O(1)=
\ln\lnc\left(1+O\left(\frac1{\ln\lnc}\right)\right).
\end{align}
Using formulas (\ref{formula poryadok lnN3}) and (\ref{formula poryadok lnlnN3}), we obtain
\begin{align*}
\frac{\ln\ln N_3}{\ln N_3}=&\frac{c\ln\lnc}{\lnc}\cdot\frac{1}{1+\R(N)}\left(1+O\left(\frac1{\ln\lnc}\right)\right) =\frac{c\ln\lnc}{\lnc}\left(1+O\left(\frac1{\lnc}\right)\right)\left(1+O\left(\frac1{\ln\lnc}\right)\right)=\\
&=\frac{c\ln \lnc }{\lnc}+O\left(\frac{1}{\ln\lnc }\right).
\end{align*}
From this and from the condition $\|c\|\ge  3\left(2^{[c]+1}-1\right)\dfrac{c\ln \lnc}{\lnc}$, it follows that
\begin{align*}
&\|c\|\ge 
3\left(2^{[c]+1}-1\right)\left(\frac{\ln \ln (N_1+H_1)}{\ln (N_1+H_1)}+O\left(\frac{1}{\lnc }\right)\right),
\end{align*}
that is, condition (\ref{uslovie1}) is satisfied. Using relations (\ref{formula poryadok H3}), (\ref{formula H=N..}), (\ref{formula poryadok N3}) and (\ref{formula poryadok lnN3}), one can establish the validity of condition (\ref{uslovie2}):
\begin{align*}
\frac{H_3}{\sqrt{2cN_3}(\ln N_3)^2}
&=\frac{\dfrac{2N^{1-\frac{1}{2c}}\lnc^2}{c(\mu_3N)^{1-\frac1c}}\left(1+O\left(\dfrac{H}{N}\right)\right)} {\left(2c(\mu_3N)^\frac1c\left(1+O\left(\dfrac HN\right)\right)\right)^\frac12\left(\dfrac\lnc c\left(1+O\left(\dfrac1\lnc\right)\right)\right)^2}=\\
&=\frac{\sqrt{2c}}{\mu_3^{1-\frac1{2c}}}\left(1+O\left(\dfrac1\lnc\right)\right)>1.
\end{align*}

Similarly, using relations (\ref{formula poryadok lnN3}), (\ref{formula poryadok H3}) and (\ref{formula poryadok N3}), we obtain:
\begin{align}
\frac{\ln^2 N_3}{H_3N_3^{c-1}}=&\frac{\left(\dfrac\lnc c\left(1+O\left(\dfrac1\lnc\right)\right)\right)^2} {\dfrac{2H}{c(\mu_3N)^{1-\frac1c}}\left(1+O\left(\dfrac{H}{N}\right)\right) \left((\mu_3N)^\frac1c\left(1+O\left(\dfrac HN\right)\right)\right)^{c-1}}=\nonumber\\
&=\frac{\lnc^2}{2cH}\left(1+O\left(\dfrac1\lnc\right)\right)=\varkappa +O\left(\frac{\varkappa}{\lnc} \right). \label{ae asiptotika}
\end{align}
From the condition $\alpha \in \m_+=[\varkappa,0.5]$, condition (\ref{uslovie3}) follows.

Thus, according to Lemma~\ref{Lemma otsenka S_c(alpha ;x,y)}, taking into account relation (\ref{formula poryadok H3}), we have
$$
S_c(\alpha;N_3,H_3)\ll\frac{H_3}{(\ln N_3)^2}\ll \frac{H}{N^{1-\frac 1c}\lnc^2}.
$$
Substituting this estimate into (\ref{Integral I(M2)}), we find
$$
I(\m_+)\ll \frac{H}{\lnc }\cdot \max_{\alpha \in \m_+}|S_c(\alpha;N_3,H_3)|\ll\frac{H^2}{N^{1-\frac 1c}\lnc^3}.
$$
The absolute values of the integrals $I(\m_+)$ and $I(\m_-)$ coincide, therefore the last estimate is also valid for $I(\m_-)$.

\subsection{Evaluation of the integral $I(\M)$}
By the definition of the integral $I(\M)$, we have:
\begin{align}\label{integral po M-1}
&I(\M )=\int_{-\varkappa}^{\varkappa} \F(\alpha)e(-\alpha N)d\alpha,\qquad \qquad \varkappa=\frac{\lnc^2}{2cH},\\ &\F(\alpha)=\S_1(\alpha;N_1,2H)\S_1(\alpha;N_2,2H)S_c(\alpha;N_3,H_3).\nonumber
\end{align}

Using the condition $c\ge\dfrac43\left(1+52\dfrac{\ln\lnc}{\lnc}\right)$ and the formula for the sum of an infinitely decreasing geometric progression whose ratio is $\dfrac{52\ln\lnc}{\lnc}$, we obtain
\begin{align}
\frac38-\frac{1}{2c}&=\frac38\left(1-\frac{4}{3c}\right)\ge\frac38\left(1-\left(1+\frac{52\ln\lnc}{\lnc}\right)^{-1}\right)= \frac38\sum_{k=1}^\infty(-1)^{k-1}\left(\frac{52\ln\lnc}{\lnc}\right)^k>\nonumber\\ &>\frac38\left(\frac{52\ln\lnc}{\lnc}-\frac{52^2\ln^2\lnc}{\lnc^2}\right)= 18.5\frac{\ln\lnc}{\lnc}+\frac{\ln\lnc}{\lnc}\left(1-1014\frac{\ln\lnc}{\lnc}\right)>18.5\frac{\ln\lnc}{\lnc}.\label{formula neravenstvo 3/8-1/2c>..}
\end{align}

To obtain an asymptotic formula for the function $\F(\alpha; N, H)$, we first determine the asymptotic behavior of the sum $S_1(\alpha; N_k, 2H)$ for $k = 1, 2$. Setting
$$
x=N_k=\mu_kN+H,\qquad y=2H,\qquad A=3,
$$
we apply Lemma~\ref{Lemma poved LinKorTrigSum v okr nulya} to these sums. We verify the following two conditions of this lemma:
\begin{equation}\label{formula uslovie M}
2H\ge (\mu_kN+H)^\frac58(\ln(\mu_kN+H))^{19.5},\qquad \varkappa=(2cH)^{-1}\,\lnc^2\le\frac{\mu_kN+H}{8\pi H^2}.
\end{equation}
Using condition (\ref{formula H=N..}) and, for brevity, introducing the notation
$$
c_2=2\left(\mu_k+\frac{H}{N}\right)^{-\frac58}\left(1+\frac{\ln\left(\mu_k+\frac{H}{N}\right)}{\lnc}\right)^{-19.5},
$$
and then applying inequality (\ref{formula neravenstvo 3/8-1/2c>..}), we successively obtain
\begin{align*}
\frac{2H}{(\mu_kN+H)^\frac58(\ln(\mu_kN+H))^{19.5}}=\frac{2N^{1-\frac{1}{2c}}\lnc^2}{(\mu_kN+H)^\frac58(\ln(\mu_kN+H))^{19.5}}= c_2N^{\frac38-\frac{1}{2c}}\lnc^{-17.5}=\\
=c_2\exp\left(\left(\frac38-\frac{1}{2c}\right)\lnc-17.5\ln\lnc\right)>c_2\exp\left(18.5\frac{\ln\lnc}{\lnc}\cdot\lnc-17.5\ln\lnc\right)=c_2\lnc,
\end{align*}
that is, the first condition in (\ref{formula uslovie M}) is satisfied.

Further, using the condition (\ref{formula H=N..}), we show that the second condition (\ref{formula uslovie M}) is also satisfied:
\begin{align*}
\frac{\varkappa}{x(2\pi y)^{-2}}=\frac{(2cH)^{-1}\lnc^2}{(\mu_k N+2H)(2\pi\cdot 2H)^{-2}}=\frac{8\pi^2H\lnc^2}{c(\mu_k N+2H)}=
\frac{8\pi^2}{c\left(\mu_k+\frac{2H}N\right)}\cdot\frac{\lnc^2}{N^\frac1{2c}}<1.
\end{align*}
Thus, both conditions of Lemma~\ref{Lemma poved LinKorTrigSum v okr nulya} are satisfied, and therefore we have
\begin{align*}
S_1(\alpha;N_k,2H )&=\frac{\sin2\pi\alpha H}{\pi\alpha}e\left(\alpha\mu_kN\right) +O\left(\frac{H}{\lnc^3}\right).
\end{align*}
Using Lemma \ref{Lemma S(alpha;Nk,H=S1(alpha;Nk+H,2H)}, the sum $\S_1(\alpha; N_k, 2H)$, $k = 1,2$, is expressed in terms of the sum $S_1(\alpha; N_k, 2H)$, and then, applying the last formula, we obtain
\begin{align*}
\S_1(\alpha;N_k,2H)&=\frac{S_1(\alpha;N_k,2H)}{\ln(\mu_kN)}+O\left(\frac{H^2}{N\lnc}\right)=\frac1{{\ln(\mu_kN)}}\frac{\sin2\pi\alpha H}{\pi\alpha}e(\alpha\mu_kN)+O\left(\frac{H}{\lnc^4}\right).
\end{align*}

From this, and from the relation $\dfrac{\sin2\pi\alpha H}{\pi\alpha}\ll H$, we obtain
\begin{align}
\S_1(\alpha;N_1,2H)\S_1(\alpha;N_2N,2H)=&\frac{\sin^22\pi\alpha H} {\pi^2\alpha^2}\frac{e(\alpha(\mu_1+\mu_2)N)}{\ln(\mu_1N)\ln(\mu_2N)}+O\left(\frac{H^2}{\lnc^5}\right).\label{formula asimp S32(alpha;N1,2H1)}
\end{align}
Now, using Lemma~\ref{Lemma poveda S_c(alpha ;x,y) v okr nulya} with $A = 2$, $x = N_3$, $y = H_3$, we find the asymptotic behavior of the sum
$$
S_c(\alpha;N_3,H_3)=\sum_{N_3-H_3<n\le N_3}e(\alpha [n^c]), \qquad \alpha \in [-\varkappa,\ \varkappa].
$$
Conditions (\ref{uslovie1}) and (\ref{uslovie2}) of Lemma~\ref{Lemma otsenka S_c(alpha ;x,y)}, whose validity was demonstrated in the estimation of $I(\m_+)$, coincide with the conditions of Lemma~\ref{Lemma poveda S_c(alpha ;x,y) v okr nulya}. In estimating the sum $S_c(\alpha,N_3,H_3)$ for $\alpha$ from the set $\m_+$ in (\ref{ae asiptotika}), it was shown that
$$
\frac{\ln^2 N_3}{H_3N_3^{c-1}}=\varkappa +O\left(\frac{\varkappa}{\lnc} \right).
$$
From this formula and the condition $\alpha \in \M=[-\varkappa,\varkappa]$, it directly follows that the condition
\begin{align*}
& |\alpha|\le\frac{\ln^2N_3}{H_3N_3^{c-1}}.
\end{align*}
is satisfied.

Therefore, according to Lemma~\ref{Lemma poveda S_c(alpha ;x,y) v okr nulya}, we have
\begin{align}\label{formula Sc(alpha;N3,H3)-2}
S_c(\alpha;N_3,H_3)=\frac{\sin\pi\alpha}{\pi\alpha}\int_{N_3-H_3}^{N_3}e(\alpha(t^c-0,5))dt +O\left(\frac{H_3|\sin\pi\alpha|}{(\ln(N_3)^2}\right).
\end{align}
Estimating the remainder term of the last formula using relations (\ref{formula poryadok lnN3}) and (\ref{formula poryadok H3}), we obtain
\begin{align*}
&\frac{H_3|\sin\pi\alpha|}{(\ln(N_3)^2}\le\frac{H_3\sin\pi \varkappa}{(\ln(N_3)^2}\ll\frac{H_3}{\lnc^2 }\cdot \sin \left(\frac{\pi \lnc^2}{2cH}\right)\ll
 \frac{H}{N^{1-\frac 1c}\lnc^2 }\cdot \frac{ \lnc^2}{H}\ll\frac1{N^{1-\frac 1c}}.
\end{align*}
Substituting this estimate into the right-hand side of formula (\ref{formula Sc(alpha;N3,H3)-2}), we obtain
\begin{align}
S_c(\alpha;N_3,H)&=\frac{1-e(-\alpha )}{2\pi i\alpha}\int_{N_3-H_3}^{N_3}e(\alpha t^c)dt+O\left(\frac1{N^{1-\frac 1c}}\right)=\nonumber\\
&=\frac{H_3(1-e(-\alpha))}{2\pi i\alpha}\gamma_c(\alpha;N_3,H_3)+O\left(\frac1{N^{1-\frac 1c}}\right),\label{Sc(alpha;N,H)-3}
\\
&\gamma_c(\alpha;N_3,H_3)=\int_{-0,5}^{0,5}e\left(\alpha(N_3+H_3(t-0.5))^c\right)dt.\label{formula opr gamma(alpha;N3,H3)-1}
\end{align}
We now find an asymptotic formula for the integral $\gamma_c(\alpha;N_3,H_3)$. Using the relations
\begin{align*}
&N_3^c=\mu_3N+H,\qquad \qquad H_3^2N_3^{c-2}\ll\frac{H^2}{N^{2-\frac2c}}\cdot N^\frac{c-2}c=\frac{H^2}N,\\
&cH_3N_3^{c-1}=c\left(\frac{2H}{c(\mu_3N)^{1-\frac1c}}+O\left(\frac{H^2}{N^{2-\frac1c}}\right)\right)\left(\mu_3N\right)^\frac{c-1}c \left(1+O\left(\frac{H}{N}\right)\right)=2H+O\left(\frac{H^2}{N}\right),
\end{align*}
which follow from relations (\ref{formula poryadok H3}) and (\ref{formula poryadok N3}), we obtain
\begin{align*}
(N_3+H_3(t&-0.5))^c=N_3^c\left(1+\frac{H_3(t-0.5)}{N_3}\right)^c= N_3^c\left(1+\frac{cH_3(t-0.5)}{N_3}+O\left(\frac{H_3^2}{N_3^2}\right)\right)=\\
&= N_3^c+cH_3N_3^{c-1}(t-0.5)+O\left(H_3^2N_3^{c-2}\right)=\\
&=\mu_3N+H+2H(t-0.5)+O\left(\frac{H^2}{N}\right)=\mu_3N+2Ht+\R_2,\qquad \R_2\ll\frac{H^2}N.
\end{align*}

From this and from (\ref{formula opr gamma(alpha;N3,H3)-1}), taking into account that $e(\alpha\R_2)-1\ll|\alpha|\R_2$ and $|\alpha|\ll\lnc^2H^{-1}$, we obtain
\begin{align}\label{formula opr gamma(alpha;N3,H3)-2}
&\gamma(\alpha;N_3,H_3)=\int_{-0,5}^{0,5}e\left(\alpha\left(\mu_3N+2Ht+\R_2\right)\right)dt
=e\left(\mu_3N\alpha\right)\frac{\sin\left(2\pi\alpha H\right)}{2\pi\alpha H}+O\left(\frac{H\lnc^2}N\right).
\end{align}
Using the Taylor formula for the functions $\cos2\pi\alpha$ and $\sin2\pi\alpha$ in a neighborhood of zero $|\alpha|\le \varkappa$, we find
\begin{align}
&\frac{1-e(-\alpha)}{2\pi i \alpha}=\frac{1-\cos 2\pi \alpha+i\sin 2\pi \alpha }{2\pi i \alpha}=\frac{1-(1+O(\alpha^2))+ 2\pi i\alpha +O(\alpha^3)}{2\pi i \alpha}
=1+O\left(\frac{\lnc^2}H\right). \label{formula 1-e(alpha)..}
\end{align}
Further, multiplying termwise formulas (\ref{formula opr gamma(alpha;N3,H3)-2}) and (\ref{formula 1-e(alpha)..}) and using the inequality $\left|\dfrac{\sin\left(2\pi\alpha H\right)}{2\pi\alpha H}\right|\ll 1$, we obtain
\begin{align*}
\frac{1-e(-\alpha)}{2\pi i \alpha}\gamma(\alpha;N_3,H_3)
=(e\left(\mu_3N\alpha\right)\frac{\sin\left(2\pi\alpha H\right)}{2\pi\alpha H}+O\left(\frac{H\lnc^2}{N}\right).
\end{align*}
From this and from (\ref{Sc(alpha;N,H)-3}) we obtain
\begin{align}
S_c(\alpha;N_3,H)&=H_3\frac{\sin(2\pi\alpha H)}{2\pi\alpha H}e(\mu_3N\alpha)+ O\left(\frac{H^2\lnc^2}{N^{2-\frac1c}}\right). \label{Sc(alpha;N,H)-4}
\end{align}
Multiplying termwise formulas (\ref{formula asimp S32(alpha;N1,2H1)}) and (\ref{Sc(alpha;N,H)-4}), and then estimating the remainder term of this product, denoted by $\R_4$, using the inequality $|\sin(2\pi \alpha H)| \ll |2 \pi \alpha H|$ and relation (\ref{formula poryadok H3}), we obtain
\begin{align*}
\F(\alpha)&=\frac{H_3}{2H\ln(\mu_1N)\ln(\mu_2N)}\cdot\frac{\sin^32\pi\alpha H}{\pi^3\alpha^3}e(\alpha(\mu_1+\mu_2+\mu_3)N)+\R_4,\\
\R_4&\ll\frac{\sin^22\pi\alpha H }{\pi^2\alpha^2\lnc^2}\cdot\frac{H^2\lnc^2}{N^{2-\frac1c}}+\frac{H_3|\sin(2\pi H\alpha)|}{|2\pi H\alpha|}\frac{H^2}{\lnc^5}+\frac{H^4}{N^{2-\frac1c}\lnc^3}\ll\\
&\ll\frac{H^4}{N^{2-\frac1c}}+\frac{H^3}{N^{1-\frac1c}\lnc^5}+\frac{H^4}{N^{2-\frac1c}\lnc^3}
 \ll\frac{H^3}{N^{1-\frac1c}\lnc^5}.
\end{align*}

Substituting the expression for the function $\F(\alpha)$, that is, the right-hand side of the last formula, into (\ref{integral po M-1}), we obtain
\begin{align}
I(\M )&=\frac{H_3}{2H\ln(\mu_1N)\ln(\mu_2N)}J(H)+O\left(\frac{H^2}{N^{1-\frac1c}\lnc^3}\right),\label{integral po M-2}\\
& J(H)=\int_{-\varkappa}^{\varkappa}\frac{\sin^32\pi\alpha H}{\pi^3\alpha^3}d\alpha.\nonumber
\end{align}
Replacing $J(H)$ by a close improper integral independent of $\lnc$, and using the formula (see \cite{Uitiker Vatson}, p.~174)
\begin{align*}
\int\limits_0^\infty\frac{\sin^nmu}{u^n}du=&\frac{\pi
m^{m-1}}{2^n(n-1)!}\left[n^{n-1}-\frac{n}{1!}(n-2)^{n-1}+\frac{n(n-1)}{2!}(n-4)^{n-1}
+\ldots\right],
\end{align*}
for $m=1$ and $n=3$, we find
\begin{align*}
J(H)&
=\frac{8H^2}{\pi}\int_0^{2\pi \varkappa H}  \frac{\sin^3 u }{u^3}du
=\frac{8H^2}{\pi}\left(\int_0^\infty\frac{\sin^3  }{u^3}du-\int_{2\pi \varkappa H}^\infty  \frac{\sin^3 u }{u^3}du\right)=
\\
&=\frac{8H^2}{\pi}\int_0^\infty  \frac{\sin^3 u }{u^3}du+O\left(\frac{H}{\lnc^6}\right)=3H^2+O\left(\frac{H^2}{\lnc^6}\right).
\end{align*}
Substituting the value of the integral $J(H)$ into formula (\ref{integral po M-2}), we obtain
$$
I(\M)=\frac{3HH_3}{2\ln(\mu_1N)\ln(\mu_2N)}+O\left(\frac{H^2}{N^{1-\frac1c}\lnc^A}\right)+O\left(\frac{H^2}{N^{1-\frac1c}\lnc^8}\right).
$$

Using formula (\ref{formula poryadok H3}) and the relation
\begin{align*}
\frac1{\ln(\mu_kN)}-\frac1{\lnc}&=\frac{-\ln\mu_k}{(\lnc-\ln\mu_k)\lnc}\ll\frac1{\lnc^2},
\end{align*}
we obtain
$$
I(\M )=\frac{3H^2}{c(\mu_3N)^{1-\frac1c}\lnc^2}+O\left(\frac{H^2}{N^{1-\frac1c}\lnc^3}\right).
$$
The theorem is proved.


\end{document}